\documentclass[journal,twoside,web]{ieeecolor}
\usepackage{lcsys}                         

\usepackage{hyperref}%
\hypersetup{colorlinks=true, linkcolor=blue, breaklinks=true, urlcolor=blue}
\usepackage{doi}

\usepackage[all]{xy}
\usepackage{amsmath}
\usepackage{amscd} 
\usepackage{mathrsfs}
\usepackage{amssymb}
\usepackage{tikz}
\usepackage[utf8]{inputenc}
\usepackage[yyyymmdd]{datetime}
\usepackage{graphicx}
\usepackage{subcaption}
\usepackage{lipsum}
\usepackage{multicol}
\usepackage{mathtools}
\usepackage{multirow}
\usepackage{threeparttable}
\usepackage{pdfpages}
\usepackage{multicol} 
\usepackage{mathtools}
\usepackage{version,xspace}
\usepackage{url,doi}
\usepackage{amsmath}

\newcommand{\suchthat}{\;\ifnum\currentgrouptype=16 \middle\fi|\;}
\newcommand{\scirc}{\raise1pt\hbox{$\,\scriptstyle\circ\,$}}

\newcommand{\real}{\mathbb{R}}

\DeclareMathOperator*{\argmin}{argmin}

\renewcommand{\top}{\mathsf{T}}  

\DeclareSymbolFont{bbold}{U}{bbold}{m}{n}
\DeclareSymbolFontAlphabet{\mathbbold}{bbold}
\newcommand{\vect}[1]{\mathbbold{#1}}

\newcommand\oprocendsymbol{\hbox{$\triangle$}}
\newcommand\oprocend{\relax\ifmmode\else\unskip\hfill\fi\oprocendsymbol}

\usepackage[prependcaption,colorinlistoftodos]{todonotes}

\newcommand{\Lip}{{\rm Lip}}
\newcommand{\osL}{{\rm osL}}

\usepackage{xcolor}
\usepackage{xparse}
\usepackage{mathrsfs}
\usepackage{ stmaryrd }
\usepackage{amssymb}
\usepackage{graphics} % for pdf,  bitmapped graphics files
\usepackage{amsmath} % assumes amsmath package installed
\usepackage{amssymb}  % assumes amsmath package installed
\usepackage{hyperref}%
\hypersetup{colorlinks=true,  linkcolor=blue,  breaklinks=true,  urlcolor=blue,  
citecolor=blue}

\usepackage{url}
\usepackage[prependcaption, colorinlistoftodos]{todonotes}

\renewcommand{\theenumi}{(\roman{enumi})}

\DeclareSymbolFont{bbold}{U}{bbold}{m}{n}
\DeclareSymbolFontAlphabet{\mathbbold}{bbold}
\usepackage{soul}

\usepackage{centernot}
\usepackage{cancel}[Smaller]
\NewDocumentCommand\seminorm{mg}{\lVert  #1 \bcancel{\rVert} \IfNoValueF{#2}{_{#2}}}

\NewDocumentCommand\bigseminorm{mg}{\Big\lvert\!\Big\lvert\!\Big\lvert #1 
\Big\rvert\!\Big\rvert\!\Big\rvert \IfNoValueF{#2}{_{#2}}}
\newtheorem{theorem}{Theorem}

\newtheorem{lemma}[theorem]{Lemma}
\newtheorem{definition}[theorem]{Definition}

\newtheorem{corollary}[theorem]{Corollary}
 
\newtheorem{assumption}{Assumption}
\newtheorem{problem}{Problem}

\NewDocumentCommand\norm{mg}{\lVert #1 \rVert \IfNoValueF{#2}{_{#2}}}
\newcommand{\bullobib}{\bibliographystyle{plainurl}\bibliography{alias,New,Main,FB}}
\usepackage{textcomp}

\definecolor{gnred}{RGB}{255,91,89}

\definecolor{gnred1}{RGB}{71,0,0} % 470000  darker
\definecolor{gnred2}{RGB}{117,0,0} % 750000
\definecolor{gnred3}{RGB}{164,0,0} % a40000
\definecolor{gnred4}{RGB}{211,0,0} % d30000
\definecolor{gnred5}{RGB}{255,0,0} % FF0000
\definecolor{gnred6}{RGB}{255,42,34} % FF2a22
\definecolor{gnred7}{RGB}{255,91,89} % ff5b59 --- favorite

\definecolor{gnblue1}{RGB}{0,36,71}   % 002447  darker
\definecolor{gnblue2}{RGB}{0,60,118}  % 003c76
\definecolor{gnblue3}{RGB}{0,85,164}  % 0055A4
\definecolor{gnblue4}{RGB}{0,108,212} % 006CD4
\definecolor{gnblue5}{RGB}{0,133,255}  % 0085ff
\definecolor{gnblue6}{RGB}{35,156,255} % 239cff
\definecolor{gnblue7}{RGB}{88,177,255} % 58b1ff

\definecolor{gnbrown1}{RGB}{71,27,0}  % 471b00 darker
\definecolor{gnbrown2}{RGB}{117,45,0} % 752d00
\definecolor{gnbrown3}{RGB}{164,62,0} % a43e00
\definecolor{gnbrown4}{RGB}{211,80,0} % d35000
\definecolor{gnbrown5}{RGB}{255,97,0} % ff6100
\definecolor{gnbrown6}{RGB}{255,127,26} % ff7f1a
\definecolor{gnbrown7}{RGB}{255,155,86} % ff9b56

\usepackage[noend]{algorithmic}
\usepackage{algorithm}
\algsetup{indent=2em}

\newsavebox{\ieeealgbox}

\usepackage{algorithmic}
\usepackage{relsize}

\title{\LARGE \bf Contractivity of the Method of Successive Approximations \\ for Optimal
  Control}

\author{Kevin D. Smith, \IEEEmembership{Graduate Student Member, IEEE}, and Francesco 
Bullo, \IEEEmembership{Fellow, IEEE}
\thanks{This work was in part supported by AFOSR projects FA9550-22-1-0059 and FA9550-21-1-0203.}%
\thanks{Kevin D.\ Smith and Francesco Bullo are with the Center for Control, Dynamical 
Systems, and Computation, UC Santa Barbara, Santa Barbara, CA 93101 USA. {\tt\small 
kevinsmith@ucsb.edu, bullo@ucsb.edu}.}}

\DeclareMathOperator{\MSA}{MSA}

\makeatletter
\newcommand{\oset}[3][0ex]{%
  \mathrel{\mathop{#3}\limits^{
    \vbox to#1{\kern-2\ex@
    \hbox{$\scriptstyle#2$}\vss}}}}
\makeatother
\newcommand{\LambBack}{\lambda^{\leftarrow}}

\newcommand{\rev}[1]{#1}

\begin{document}

\maketitle
\thispagestyle{empty}
\pagestyle{empty}
 
\begin{abstract}
  Strongly contracting dynamical systems have numerous properties (e.g.,
  incremental ISS), find widespread applications (e.g., in controls and
  learning), and their study is receiving increasing attention. This work starts with the 
  simple observation that, given a strongly contracting system, its adjoint dynamical 
  system is also strongly contracting, with the same rate, with respect to the dual norm, 
  under time reversal. As main implication of this dual contractivity, we show that the 
  classic Method of Successive Approximations (MSA), an indirect method in optimal 
  control, is a contraction mapping for short optimization intervals or large contraction 
  rates. Consequently, we establish new convergence conditions for the MSA algorithm, 
  which further imply uniqueness of the optimal control and sufficiency of Pontryagin's 
  minimum principle under additional assumptions.
  
%  Our approach therefore provides a constructive
%  approach to the existence and uniqueness of solutions to the Pontryagin's
%  minimum principle.
\end{abstract} 

%  \show{abstract}

\section{Introduction} 

Optimal control is generally a difficult problem, and with the exception of some 
analytically tractable cases, it must be solved numerically. Numerical approaches broadly 
fall into two categories: direct and indirect methods. Direct methods, like direct 
collocation and direct shooting methods \cite{OVS:93, AVR:09, JTB:10}, discretize and 
approximate the state and/or control to encode the problem as a nonlinear program. Due to 
their relative simplicity, robustness, and the wide availability of software 
implementations, direct methods tend to be favored in modern times \cite[\S 
4.3]{JTB:10}, \cite{BAC:12}.  

Indirect methods are an older class of methods based on Pontryagin's minimum principle 
(PMP), which gives a necessary condition for optimality of a control signal. PMP states 
that the optimal trajectory must solve a two-point boundary problem, together with a 
costate, and that the optimal control minimizes a Hamiltonian function at each 
point in time. Indirect methods search for an input, state trajectory, and costate 
trajectory that satisfy PMP. Many direct methods, including shooting and collocation, can 
also be applied as indirect methods to the PMP boundary value problem \cite{HBK:18}. 
Another approach is the Method of Successive Approximations (MSA) \cite{FLC-AAL:82}, also 
called the Forward-Backward-Sweep algorithm \cite{SL-JTW:07}, which is the main topic of 
this letter. 
 
MSA \cite{HJK-REK-HGM:61, IAK-FLC:63, VVA:68} and its variants \cite{SKM:66, FLC-AAL:82} 
are classic approaches that have received renewed attention in the machine learning 
community \cite{QL-SH:18, QL-LC-CT-WE:18, LB-NAF-TA:22} as alternatives to gradient 
descent for 
training residual neural networks (ResNets). Indeed, a new thrust of machine learning 
research is to apply control-theoretic
techniques to the training of ResNets by viewing these models as forward Euler
discretizations of continuous-time control systems \cite{WE:17, 
DZ-TZ-YL-ZZ-BD:19,PT-BG:22}. Within this framework, training the ResNet can be viewed as 
an optimal control problem. As argued in \cite{QL-SH:18, QL-LC-CT-WE:18}, MSA (and its 
variants) allow for error and convergence analysis and can lead to better 
training dynamics than gradient descent. 

Unfortunately, MSA does not always converge, a problem that is still the subject of 
ongoing research. In \cite{MMcA-LM-WH:12}, the authors prove convergence criteria based 
on boundedness and Lipschitz assumptions. Similar bounds are established in 
\cite{QL-SH:18, QL-LC-CT-WE:18}. This letter provides a new set of convergence criteria 
\rev{when MSA is applied to strongly contracting dynamical systems.}

The contributions of this letter are as follows. First, in \S\ref{sect:adjoint}, we study
the adjoints of nonlinear systems that arise in optimal control theory. We show that
adjoints of contracting systems under time reversal are also contracting with the same
rate, albeit with respect to the dual norm. This 
property allows us to prove Gr\"onwall-like and ISS-like bounds on the adjoint dynamics. 
\S\ref{sect:msa} applies 
these bounds to analyze MSA. Assuming Lipschitz continuity of all relevant maps in the 
optimal control problem, we obtain a bound on the Lipschitz constant of each MSA 
iteration. This Lipschitz constant becomes arbitrarily small in the limits of short 
optimization intervals and large contraction rates, thereby establishing conditions for 
when the iteration is a contraction mapping. With an additional assumption of pointwise 
uniqueness of the minimizer of the Hamiltonian, we show that these conditions also lead 
to uniqueness of the optimal control and sufficiency of PMP. Finally, in \S 
\ref{sect:example}, we provide an illustrative example.

\section{Preliminaries}
 
\subsection{Contracting Dynamics over Normed Vector Spaces}

Let $\norm{\cdot}: \real^n \to \real_{\ge 0}$ be a norm. The dual norm
$\norm{\cdot}_\star: \real^n \to \real_{\ge 0}$ is the norm $\norm{x}_\star =
\sup_{\norm{y} \le 1} y^\top x$.  Given a matrix $A \in \real^{n \times n}$, the induced
norm of $A$ is $\norm{A} = \sup_{\norm{x} = 1} \norm{Ax}$ and the induced logarithmic norm
of $A$ is
\[
	\mu(A) = \lim_{\rev{\alpha} \to 0^+} \frac{\norm{I_n + \rev{\alpha} A} - 
	1}{\rev{\alpha}}.
\]
Explicit formulas for the induced (logarithmic) norms are known for the standard $p \in
\{1, 2, \infty\}$ norms on $\real^n$ \cite[\S 2.4]{FB:22-CTDS}.

A map $T: X \to Y$ between normed spaces $(X, \norm{\cdot}_X)$ and $(Y, \norm{\cdot}_Y)$
is Lipschitz continuous if a constant $\ell \ge 0$ exists such that $\norm{T(x) - T(\bar 
x)}_Y \le \ell \norm{x - \bar x}_X$ for all $x, \bar x \in X$. The \emph{\rev{minimal} 
Lipschitz constant}
$\Lip(T)$ is the infimum over $\ell$ that satisfy this inequality. If $T$ is continuously
differentiable, then $\Lip(T) = \sup_{x \in X} \norm{D_x \rev{T}(x)}$, \rev{where $D_x 
T(x)$ denotes the Jacobian matrix of $T$.} Furthermore, if 
\rev{$X = Y = \real^n$}, then the \emph{one-sided 
Lipschitz constant} of $T$ is $\osL(T) = \sup_{x \in X} \mu(D_x 
\rev{T} (x))$. A dynamical system $\dot x = f(t, x, 
\dots)$ \rev{with a continuously differentiable vector field $f: \real^n \to \real^n$} is 
said to be \emph{strongly infinitesimally contracting} with rate $c > 0$ if the map $x 
\mapsto f(t, x, \dots)$ is \rev{uniformly} one-sided Lipschitz with constant $-c$ for all 
\rev{$t$} and for all inputs. 

Strongly contracting systems enjoy numerous properties. As a useful example, we state the 
following lemma without proof (as it slightly generalizes \cite[Theorem 3.15, Corollary 
3.16]{FB:22-CTDS}).

\begin{lemma}[Gr\"onwall comparison lemma] \label{lem:iss}
  Consider a dynamical system
  \begin{equation} \label{eq:f}
    \dot x(t) = f(t, x(t), u_1(t), \dots, u_m(t)), \qquad \forall t \ge 0,
  \end{equation}
  with $x(t) \in \real^n$ and inputs $u_i \in U_i \subseteq \real^{k_i}$ for $i 
  \in \{1, 2, \dots, m\}$. Let $\norm{\cdot}$ be a norm on $\real^n$, and let 
  $\norm{\cdot}_{U_i}$ be norms on $U_i$. Assume that
  \begin{enumerate}
  \item \label{item:gron-contract} the system~\eqref{eq:f} is strongly infinitesimally 
  contracting with rate $c > 0$, and
  \item \rev{for each $i \in \{1, 2, \dots, m\}$, the maps $u_i \mapsto f(t, x, u_1, 
  \dots, u_i, \dots, u_m)$ \rev{are uniformly} Lipschitz continuous with constant 
  $\ell_{f, U_i}$ for all $t \ge 0$, $x \in \real^n$, and $u_j \in U_j$ with $j \ne i$.}
  \end{enumerate}
  Let $(u_1, \dots, u_m)$ and $(\bar u_1, \dots, \bar u_m)$ be input signals, and 
  let $x, \bar x$ be the corresponding trajectories of \eqref{eq:f}. For all $t \ge 0$, 
  \begin{align} \label{eq:gron-bound}
    \begin{split}
      &\norm{x(t) - \bar x(t)} \le e^{-c t} \norm{x(0) - \bar x(0)} \\
      &\qquad + \sum_{i = 1}^m \ell_{f, \rev{U}_i} \int_0^t e^{-c (t - \tau)} 
      \norm{u_i(\tau) - \bar u_i(\tau)}_{U_i} d\tau .
    \end{split}
  \end{align}  
\end{lemma} 

\rev{Note that \eqref{eq:gron-bound} still holds when $c \le 0$, i.e., for 
expansive systems with a bounded rate of expansion; however, we do not consider such 
systems in this letter.}

\smallskip

\subsection{Optimal Control}

We study the following optimal control problem:

\rev{
\begin{problem}[Optimal control problem] \label{prob:opt}
	Consider a dynamical system 
	\begin{equation} \label{eq:sys}
		\dot x(t) = f(t, x(t), u(t)), \quad x(0)=x_0 \in \real^n, 
	\end{equation}
	where $f$ is continuous in all arguments and continuously differentiable in the 
	second and third arguments. Further consider a cost functional 
	\begin{equation} \label{eq:cost}
			J[u] = \int_0^{T} \phi(t, x(t), u(t)) \, dt + \psi(x(T)),
	\end{equation}
	where $\phi: [0, T] \times \real^n \times U \to \real$ is a \emph{running cost} that 
	is differentiable in the second argument, and $\psi: \real^n \to \real$ is a 
	differentiable \emph{terminal cost}. Let $\mathcal U = \left\{ u: [0, T] \to U 
	~\text{s.t.}~ u~\text{measurable} \right\}$ be a space of permissible control 
	signals, where $T > 0$ and $U \subseteq \real^k$ is a compact set containing $\vect 
	0_k$. The \emph{optimal control problem} is to find $u^* \in \mathcal U$ that 
	minimizes $J[u^*]$.
\end{problem} 
}
\smallskip

%More complex optimal control problems have been studied, including constraints on the 
%terminal state and free time horizons, but this letter focuses on Problem 
%\ref{prob:opt}. 
An elementary necessary condition for the optimality of a control \rev{is} Pontryagin's 
minimum principle (PMP) \cite[Theorem 5.10, Theorem 5.11]{MA-PLF:66} \cite[Theorem 
6.3.1, Theorem 6.5.1]{AB-BP:07}:
 
\begin{theorem}[Pontryagin's minimum principle]
	Let $u^* \in \mathcal U$ be an optimal control for Problem \ref{prob:opt} (if one 
	exists), and let 
	$x: [0, T] \to \real^n$ be the corresponding trajectory of \eqref{eq:sys}. For all $t 
	\in [0, T]$,
	\begin{equation} \label{eq:optimality}
		\rev{u^*}(t) \in \argmin_{\tilde u \in U} H(t, x(t), \lambda(t), \tilde u),
	\end{equation}
	where $H: \real \times \real^n \times \real^n \times U \to \real$ is the 
	\textit{Hamiltonian} 
	\begin{equation} 
		H(t, x, \lambda, u) = \lambda^\top f(t, x, u) + \phi(t, x, u)
	\end{equation} 
	and $\lambda: [0, T] \to \real^n$ is the \textit{costate} trajectory 
	\begin{equation} \label{eq:costate}
		\dot \lambda(t) = -D_x f(t, x(t), u(t))^\top \lambda(t) - \phi_x(t, 
		x(t), u(t))
	\end{equation}
	with the boundary condition $\lambda(T) = \psi_x(x(T))$.
\end{theorem} 

%\smallskip

%If $\rev{\nu} = 0$, the problem is said to be \textit{singular} or 
%\rev{\textit{abnormal}}, in which case \rev{PMP is trivially satisfied with $\lambda(t) 
%= 
%\vect 0_n$ and $H = 0$} \cite[\S 6]{MA-PLF:66}. Otherwise, we adopt the standard 
%assumption (without loss of generality) that $\rev{\nu} = 1$. Singular problems are 
%notoriously difficult, and we do not consider them in this letter.

%The existence of an optimal control can be established under standard boundedness and 
%convexity assumptions. 
%
%\begin{theorem}[Existence of optimal control] \label{thm:existence}
%	If there exists a bounded set $X \subset \real^n$ such that $x(t) \in X$ for all $t 
%	\in [0, T]$ and $u \in \mathcal U$, and for each $t \in [0, T]$ and $x \in X$, 
%	the set
%	\[ 
%		F(t, x) = \left\{
%			(\phi(t, x, u), f(t, x, u)) \rev{: u \in U}
%		\right\} \rev{\subseteq \real^{n + 1}}
%	\]
%	is convex, then an optimal control exists.
%\end{theorem}
%\smallskip
%
%\begin{proof}
%	We use an argument from 
%	\cite[\S 5.1]{AB-BP:07}: first, transform Problem \ref{prob:opt} into a Mayer problem 
%	by defining auxiliary dynamics $\dot y(t) = \phi(t, x(t), u(t))$ with $y(0) = 0$, so 
%	that $J[u] = y(T) + \psi(x(T))$. Then \cite[Theorem 
%	3.5.1]{AB-BP:07} implies that the set of reachable $(x(T), y(T))$ is compact. Hence a 
%	terminal state $(x(T), y(T))$ exists which minimizes $J$, and any $u^*$ which drives 
%	the system to this state is optimal.
%\end{proof}

\subsection{Method of Successive Approximations}

The Method of Successive Approximations (MSA) \cite{FLC-AAL:82}, also called the 
Forward-Backward Sweep algorithm \cite{SL-JTW:07}, is a basic approach to computing an 
input that satisfies PMP. The method iteratively solves the PMP 
two-point boundary value problem, then updates the control to minimize the new 
Hamiltonian at each time, as outlined in Algorithm \ref{alg:msa}.

\begin{algorithm}
	\caption{Method of Successive Approximations}
	\label{alg:msa}
	\begin{algorithmic}[1]
		\REQUIRE initial guess $u^{(0)} \in \mathcal U$
		\FOR{$i = 1, 2, \dots, N$}
			\STATE $x^{(i)} \leftarrow $ trajectory of \eqref{eq:sys} from $x_0$ with 
			input $u^{(i - 1)}$
			\STATE $\lambda^{(i)} \leftarrow $ trajectory of \eqref{eq:costate} from 
			$\lambda(T) = \psi_x(\rev{x^{(i)}}(T))$ with inputs $x^{(i)}$ and $u^{(i - 
			1)}$
			\STATE $u^{(i)}(t) \leftarrow \argmin_{\tilde u \in U} H(t, x^{(i)}(t), 
			\lambda^{(i)}(t), \tilde u)$ for all $t \in [0, T]$, ties broken arbitrarily
		\ENDFOR
		\RETURN $u^{(N)}$
	\end{algorithmic}
\end{algorithm}

The algorithm can run for a fixed number of iterations; alternatively, it may terminate 
when the difference between successive iterates $u^{(i - 1)}$, $u^{(i)}$ is within a 
specified tolerance. Note that each iteration of the algorithm maps a control $u^{(i - 
1)}$ to a new control $u^{(i)}$, so that each iteration can be thought of as an operator 
$\MSA: \mathcal U \to \mathcal U$.

\begin{definition}[MSA Operator] \label{def:msa}
	Given a control $u \in \mathcal U$, let $x: [0, T] \to \real^n$ be the 
	\rev{corresponding} trajectory of 
	\eqref{eq:sys}, and let $\lambda: [0, T] \to \real^n$ be the trajectory of 
	\eqref{eq:costate} from $\lambda(T) = \psi_x(x(T))$. Then $\MSA(u)$ is the
	control that satisfies \eqref{eq:optimality} \rev{with respect to $x(t)$ and 
	$\lambda(t)$ for all $t \in [0, T]$, with ties broken in an arbitrary deterministic 
	manner}.
\end{definition}   

\rev{Definition \ref{def:msa} is well-posed if the signal of Hamiltonian-minimizing 
controls from \eqref{eq:optimality} is measurable. When we analyze the MSA algorithm in 
\S \ref{sect:msa}, we will impose Lipschitz continuity assumptions that forbid any edge 
cases where $\MSA(u)$ is not measurable.}

\subsection{Adjoints}
\label{sect:adj}

Adjoints are familiar from linear systems theory. Given input and output Hilbert spaces 
$\mathcal S_{\rm in}, \mathcal S_{\rm out}$ and a linear system $G: \mathcal S_{\rm in} 
\to \mathcal S_{\rm out}$, the adjoint of $G$ is the unique linear system $\tilde G: 
\mathcal S_{\rm out} \to \mathcal S_{\rm in}$ such that $\langle G u, y \rangle_{\mathcal 
S_{\rm out}} = \langle u, \tilde G y \rangle_{\mathcal S_{\rm in}}$ for all $u \in 
\mathcal S_{\rm in}$ and $y \in \mathcal S_{\rm out}$. For an LTV system with the usual 
$(A, B, C, D)$ representation, the adjoint dynamics are
\begin{subequations} \label{eq:adjoint-linear}
\begin{align}
	\dot \lambda(t) &= -A(t)^\top \lambda(t) - C(t)^\top v(t) \label{eq:adjoint-state} \\
	z(t) &= B(t)^\top \lambda(t) + D(t)^\top v(t)
\end{align}
\end{subequations}
with $v \in \mathcal S_{\rm out}$ and $z \in \mathcal S_{\rm in}$ 
\cite[3.2.4]{MG-DJNL:95}. The theory of adjoints leads to the duality of
controllability and observability and of linear quadratic regulators and estimators 
\cite{OK-DSB:20}.

\section{Contractivity of the Adjoint}
\label{sect:adjoint}

This section examines adjoints of nonlinear systems. We first explain how the notion of
``adjoint'' frequently used in the optimal control literature relates to the adjoint from
linear systems. We then prove a simple yet powerful result: that the adjoint of a strongly
infinitesimally contracting system is itself strongly infinitesimally contracting, with
respect to the dual norm, when integrated backwards in time. This dual contractivity
property leads to useful bounds for the evolution of costates, to later be employed in \S 
\ref{sect:msa}.

\subsection{Adjoints of Nonlinear Systems}
Nonlinear systems do not properly have adjoints according to the definition in
\S\ref{sect:adj}. Instead, the adjoint of the system's linearized variational dynamics is 
often referred to as its adjoint \cite{HJK-REK-HGM:61, PEC-FLL-AJvdS:95}. Consider the 
nonlinear system \eqref{eq:sys} with output $y(t) = x(t)$. Let $u(t)$ be an input signal 
corresponding to a nominal trajectory $x(t)$, let $\tilde x(t)$ be the trajectory from 
$\tilde u(t)$. Linearizing the dynamics of $\delta x(t) = \tilde x(t) - x(t)$ from 
$\delta u(t) = \tilde u(t) - u(t)$,
\begin{align*}
	\dot{(\delta x)}(t) &= D_x f(t, x(t), u(t)) \delta x(t) + D_u f(t, x(t), u(t)) 
	\delta u(t) \\
	{(\delta y)(t)} &= \delta x(t)
\end{align*} 
so by \eqref{eq:adjoint-state}, the adjoint dynamics are
\begin{subequations} \label{eq:adjoint}
\begin{align} 
	\dot \lambda(t) &= -D_x f(t, x(t), u(t))^\top \lambda(t) - v(t) \label{eq:lamb} \\
	z(t) &= D_u f(t, x(t), u(t))^\top \lambda(t) \label{eq:z}
\end{align} 
\end{subequations}
where $v(t) \in V \subseteq \real^n$. Not coincidentally, the costate dynamics 
\eqref{eq:costate} from PMP are of the form
\eqref{eq:lamb}, with a forcing term $v(t) = \phi_x(t, x(t), u(t))$ from the 
running cost. Indeed, PMP can be derived from the variational linearization described 
above; see \cite[Theorem 2.3.1, Theorem 6.1.1]{AB-BP:07}.

\subsection{Contractivity of the Adjoint}

We now examine the adjoints of strongly contracting systems. When the original system is 
contracting with respect to a norm $\norm{\cdot}$, it is natural to study the adjoint 
system using the dual norm $\norm{\cdot}{\star}$, as the following lemma suggests.

\begin{lemma}[Dual Lipschitz constants] \label{lem:dual-lip} 
	Let $\norm{\cdot}: \real^n \to \real_{\ge 0}$ be a norm, and let $\norm{\cdot}_\star$ be 
	its dual norm. Let $f, g: \real^n \to \real^n$ be a pair of continuously 
	differentiable vector fields, such that $D_x f(x) = 
	D_x g(x)^\top$ for all $x \in \real^n$. Then
	\begin{enumerate}
		\item $\Lip_{\norm{\cdot}}(f) = \Lip_{\norm{\cdot}_\star}(g)$, and
		\item $\osL_{\norm{\cdot}}(f) = \osL_{\norm{\cdot}_\star}(g)$.
	\end{enumerate}
\end{lemma}   
\smallskip 

The following is an immediate consequence of Lemma \ref{lem:dual-lip}:

\begin{theorem}[\rev{Dual contraction}] \label{thm:dual}
	\rev{Consider the pair of dynamical systems \eqref{eq:sys} and \eqref{eq:lamb}.}
	Let $T > 0$ and $c > 0$, let $\norm{\cdot}$ be a norm on $\real^n$, and let 
	$\LambBack(t) = \lambda(T - t)$ be the time-reversed \rev{trajectory of 
	\eqref{eq:lamb} (where 
	we study the time-reversed dynamics due to the minus sign in the vector field)}. The 
	following are equivalent:
	\begin{enumerate}
		\item the $x(t)$ system is strongly infinitesimally contracting with respect to 
		$\norm{\cdot}$ with rate $c$, and
		\item the $\LambBack(t)$ system is strongly infinitesimally contracting
                  with respect to $\norm{\cdot}_\star$ with rate $c$.
	\end{enumerate}
\end{theorem}

\rev{
\begin{proof}
	Let $g$ be the function
	\[
		g(t, \tilde \lambda, \tilde x, \tilde u) \triangleq \frac{d \LambBack(t)}{dt} =  
		D_x f(t, \tilde x, \tilde u)^\top \tilde \lambda + v(t),
	\]
	For all fixed $t$, $\tilde \lambda$, $\tilde x$, and 
	$\tilde u$, we have $D_{\lambda} g(t, \tilde \lambda, \tilde x, \tilde u) = D_x f(t, 
	\tilde x, \tilde u)^\top$. Hence, applying Lemma \ref{lem:dual-lip} to the maps 
	$\tilde g(\lambda) = g(t, \lambda, \tilde x, \tilde u)$ and $\tilde f(x) = f(t, x, 
	\tilde u)$, we obtain $\osL_{\norm{\cdot}}(\tilde f) = 
	\osL_{\norm{\cdot}{\star}}(\tilde g)$. Thus the maps $\lambda = g(t, \lambda, \tilde 
	x, \tilde u)$ are uniformly one-sided Lipschitz with constant $-c$ with respect to 
	$\norm{\cdot}{\star}$, if and only if the maps $x \mapsto f(t, x, \tilde u)$ have the 
	same property with respect to $\norm{\cdot}$.
\end{proof}
}

\subsection{Bounds on Adjoint Dynamics} 

Theorem \ref{thm:dual} establishes that $\LambBack(t)$ is strongly contracting so long as 
the original system is strongly contracting, so we can exploit standard bounds on 
contracting systems to bound the evolution of $\lambda(t)$. Before stating these bounds, 
we impose the following two assumptions:

\begin{assumption}[Strong contractivity] \label{ass:contraction}
  The system \eqref{eq:sys} is strongly infinitesimally contracting with
  rate $c > 0$, i.e., $\osL(f(t, \tilde x, \tilde u)) \le -c$ for all $t
  \in [0, T]$ and $\tilde u \in U$. \rev{Furthermore, the trajectory of \eqref{eq:sys} on 
  the interval $[0, T]$ with $u(t) = \vect 0_k$ is bounded.} 
\end{assumption}  

\smallskip

\begin{assumption}[Lipschitz continuity, Pt. I] \label{ass:lip-1}
	For all fixed $t \in [0, T]$ and $x \in \real^n$, the map $u \mapsto f(t, \tilde 
	x, u)$ from $(\real^k, \norm{\cdot}_U)$ into $(\real^n, \norm{\cdot})$ is Lipschitz 
	with constant $\ell_{f, u}$.
\end{assumption}
\smallskip 

\rev{Strong contractivity is a fairly strong assumption. For example, if some $x^* \in 
\real^n$ is an equilibrium point of the unforced system for all $t$, strong contractivity 
implies that $x^*$ is globally exponentially stable (due to Lemma \ref{lem:iss}).} Due to 
Theorem \ref{thm:dual}, \rev{the assumption also} implies that the adjoint 
dynamics are also strongly contracting. Consequently, we can prove that all state and 
costate trajectories remain bounded.

\begin{lemma}[Boundedness of state and costate] \label{lem:bounded}
	Consider the system \eqref{eq:sys} and its adjoint \eqref{eq:costate}. If the input 
	spaces $U \subset (\real^k, \norm{\cdot}{U})$ and $V \subset (\real^n, 
	\norm{\cdot}{\star})$ are bounded, then under Assumptions 
	\ref{ass:contraction} and \ref{ass:lip-1}, there exist 
	bounded sets $X \subset (\real^n, \norm{\cdot})$ and $\Lambda \subset (\real^n, 
	\norm{\cdot}{\star})$ such that $x(t) \in X$ and $\lambda(t) \in \Lambda$ for all $t 
	\in [0, T]$ and measurable $u: [0, T] \to U$ and $v: [0, T] \to V$.
\end{lemma}
\smallskip

In the remainder of this letter, we will let $X, \Lambda \subset \real^n$ be the bounded 
sets guaranteed by Lemma \ref{lem:bounded}. In particular, the boundedness of 
$\lambda(t)$ allows us to impose additional Lipschitz continuity assumptions:

\begin{assumption}[Lipschitz continuity, Pt. II] \label{ass:lip-2}
	For all fixed $t \in [0,T]$, $\tilde x \in X$, $\tilde u \in U$, and $\tilde \lambda 
	\in \Lambda$, 
	\begin{enumerate}
		\item the map $x \mapsto D_x f(t, x, \tilde u)^\top \tilde \lambda$ from
		$(\real^n, \norm{\cdot})$ into $(\real^n, \norm{\cdot}_\star)$ is Lipshitz with 
		constant $\ell_{f_x, x}$, and
		\item the map $u \mapsto D_x f(t, \tilde x, u)^\top \tilde \lambda$ from
		$(\real^k, \norm{\cdot}_U)$ into $(\real^n, \norm{\cdot}_\star)$ is Lipschitz 
		with 
		constant $\ell_{f_x, u}$.
	\end{enumerate}
\end{assumption}  
\smallskip

We are now ready to state the first bound on the evolution of the adjoint trajectories.

\begin{theorem}[Gr\"onwall comparison of costates] \label{thm:costate-iss}
  Consider the system \eqref{eq:sys} and its adjoint \rev{\eqref{eq:lamb}} with 
  Assumptions \ref{ass:contraction}--\ref{ass:lip-2}.  Let $u, \bar u: [0, T] \to U$ and 
  $v, \bar v: [0, T] \to V$ be two pairs of measurable input signals, and let $\lambda, 
  \bar \lambda: [0, T] \to \real^n$ be the corresponding adjoint trajectories. Then for 
  all $t \ge 0$,
  \begin{align}
    \begin{split} \label{eq:costate-iss}
      &\norm{\lambda(t) - \bar \lambda(t)}_\star \le e^{-c(T - t)}
      \norm{\lambda(T) - \bar \lambda(T)}_\star \\ 
      & \;\;\;\; + \int_t^T e^{-c(\tau - t)} \norm{v(\tau) - \bar v(\tau)}_\star \, d\tau 
      \\  
      & \;\;\;\; + \ell_{f_x, u} \int_t^T e^{-c(\tau - t)} \norm{u(\tau) - \bar u(\tau)}_U \,d\tau \\
      & \;\;\; + \frac{\ell_{f_x, x} \ell_{f, u} \sinh(c(T \!-\! t))}{c} 
      \!\!\;\int_0^t\!\!\! e^{-c(T - \tau)} \norm{u(\tau) \!-\! \bar u(\tau)}_U \,d\tau \\
      & \;\;\; + \frac{\ell_{f_x, x} \ell_{f, u} e^{-c(T - t)}}{c} 
      \!\! \int_t^T\!\!\!\! \sinh(c(T \!-\! \tau)) \norm{u(\tau) \!-\! \bar 
      u(\tau)}_U\,d\tau.
    \end{split}
  \end{align} 
\end{theorem}  
\smallskip

Theorem \ref{thm:costate-iss} provides a somewhat unwieldy bound. We can sacrifice its 
sharpness to obtain a much simpler incremental ISS property.  

\begin{corollary}[Incremental ISS of adjoint systems] \label{cor:sup-bound}
  Under the same hypotheses as Theorem \ref{thm:costate-iss},
  \begin{align}
    \begin{split} \label{eq:costate-iss-sup}
      &\sup_{t \in [0, T]} \norm{\lambda(t) - \bar \lambda(t)}_\star \le \norm{\lambda(T) 
      - \bar \lambda(T)}_\star \\ 
      &\qquad \;\; + \kappa \sup_{t \in [0, T]} \norm{v(t) - \bar v(t)}_\star \\ 
      &\qquad \;\; + \left( \ell_{f_x, u} \kappa + \ell_{f_x, x} \ell_{f, u} 
      \kappa^2\right) \sup_{t \in [0, T]} \norm{u(t) - \bar u(t)}_U  .
    \end{split}
  \end{align} 
  where 
	\begin{equation} \label{eq:kappa}
	   	\kappa = c^{-1}(1 - e^{-c T}).
	\end{equation} 
\end{corollary}   

\section{Applications to Optimal Control} 
\label{sect:msa}

Here we show how the contractivity of the adjoint system leads to the contractivity of the
MSA iteration, under additional Lipschitz continuity assumptions.
 
\begin{assumption}[Lipschitz continuity of cost \rev{gradients}] \label{ass:input-lip} 
	For all fixed $t \in [0, T]$, $\tilde x \in X$, and $\tilde u \in U$, 
	\begin{enumerate}
		\item the map $x \mapsto \phi_x(t, x, \tilde u)$ from $(\real^n, \norm{\cdot})$ 
		into $(\real^n, \norm{\cdot}_\star)$ is Lipschitz with constant $\ell_{\phi_x, 
		x}$,
		\item the map $u \mapsto \phi_x(t, \tilde x, u)$ from $(\real^k, \norm{\cdot}_U)$ 
		into $(\real^n, \norm{\cdot}_\star)$ is Lipschitz with constant $\ell_{\phi_x, u}$, 
		and
		\item the map $x \mapsto \psi_x(x)$ from $(\real^n, \norm{\cdot})$ into $(\real^n, 
		\norm{\cdot}_\star)$ is Lipschitz with constant $\ell_{\psi_x, x}$.  
	\end{enumerate} 
\end{assumption} 
\smallskip

\begin{assumption}[Lipschitz continuity of the optimum] \label{ass:opt-lip}
	There exists a \rev{continuous} map $h: [0,T] \times X \times \Lambda \to U$ such that
	\begin{equation}
		h(t, x, \lambda) \in \argmin_{u \in U} H(t, x, \lambda, u) 
	\end{equation}
	for all $t \in [0, T]$, $x \in X$, and $\lambda \in \Lambda$, with ties broken in an 
	identical manner as the MSA operator, where for all fixed $t 
	\in [0, T]$, $\tilde x \in X$, and $\tilde \lambda \in \Lambda$, 
	\begin{enumerate}
		\item the map $x \mapsto h(t, x, \tilde \lambda)$ from $(\real^n, \norm{\cdot})$ 
		into $(\real^k, \norm{\cdot}_U)$ is Lipschitz with constant $\ell_{h, x}$, and
		\item the map $\lambda \mapsto h(t, \tilde x, \lambda)$ from $(\real^n, 
		\norm{\cdot}_\star)$ into $(\real^k, \norm{\cdot}_U)$ is Lipschitz with constant 
		$\ell_{h, \lambda}$.  
	\end{enumerate}
\end{assumption}  
\smallskip 
\rev{Notice that Assumption \ref{ass:opt-lip} implies that $\MSA(u)$ is measurable for 
any $u \in \mathcal U$.} With these Lipshitz assumptions, we can finally bound the 
Lipschitz constant of the MSA operator.

\begin{theorem}[Contractivity of MSA] \label{thm:main}
	Suppose that Problem \ref{prob:opt} is nonsingular and satisfies Assumptions 
	\ref{ass:contraction}--\ref{ass:opt-lip}, and consider the norm   
	$\norm{\cdot}_\mathcal U: \mathcal U \to \real_{\ge 0}$ given by 
	\begin{equation}
		\norm{u}_\mathcal U = \sup_{t \in [0, T]} \norm{u(t)}_U .
	\end{equation}
	The following are true:
	\begin{enumerate}
		\item \label{item:lip} The Lipschitz constant of an MSA iteration with respect to 
		the $\norm{\cdot}_{\mathcal U}$ norm is bounded by
		\begin{equation}
			\Lip(\MSA) \le b_1 \kappa + b_2 \kappa^2
		\end{equation} 
		where
		\begin{subequations}
		\begin{align} 
			\label{eq:b1} b_1 &= \ell_{h, x} \ell_{f, u} + \ell_{h, \lambda} \left(
				\ell_{\psi_x, x} \ell_{f, u} + \ell_{\phi_x, u} + \ell_{f_x, u}
			\right) \\
			\label{eq:b2} b_2 &= \ell_{h, \lambda} \ell_{f, u} \left(\ell_{\phi_x, x} + 
			\ell_{f_x, x} \right)
		\end{align}
		\end{subequations} 
		\item \label{item:banach} If $b_1 \kappa + b_2 \kappa^2 < 1$, then the $\MSA$ 
		operator is a 
		contraction; hence it has a unique fixed point $\hat u \in \mathcal U$, the 
		$\MSA$ iterates $u^{(i)} = \MSA^i(u^{(0)})$ converge to $\hat u$ from any initial 
		guess $u^{(0)} \in \mathcal U$, and
		\[
			\norm{u^{(i)}(t) - \rev{\hat u}(t)}_U \le \left( \frac{(b_1 \kappa + b_2 
			\kappa^2)^i}{1 - b_1 \kappa - b_2 \kappa^2} \right) \norm{u^{(1)} - 
			u^{(0)}}_{\mathcal U}
		\] 
		for all $t \in [0, T]$.
	\end{enumerate}
\end{theorem}  
\smallskip 

\begin{corollary}[Uniqueness and Sufficiency] \label{corollary}
	Under the same hypotheses as Theorem \ref{thm:main}, if additionally
	\begin{enumerate}
		\item the Hamiltonian has a unique minimizer for all $t \in [0, T]$, $x \in X$, 
		and $\lambda \in \Lambda$,
		\item an optimal control $u^*$ exists, and
		\item the time horizon $T$ is sufficiently small or the contraction rate $c$ is 
		sufficiently large that $\Lip(\MSA) < 1$,
	\end{enumerate}
	then $u^*$ is the unique optimal control, and PMP is a sufficient condition for 
	optimality. 
\end{corollary} 

\section{Example} \label{sect:example}

\begin{figure}
	\centering
	\includegraphics[width=2in]{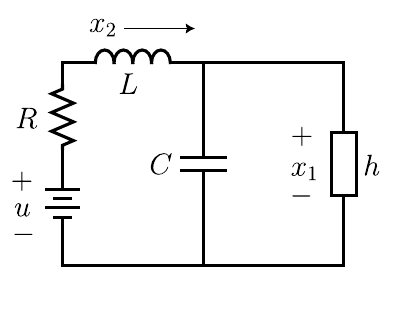}
	\caption{Diagram of the nonlinear circuit studied in Section \ref{sect:example}.}
	\label{fig:circuit}
\end{figure}

For an illustrative example the results, consider the circuit from \cite[\S 
1.2.2]{HKK:02}, depicted in Figure \ref{fig:circuit}. The circuit contains a nonlinear 
resistive element, with a current-voltage relationship $i_h = r(v_h)$ for some 
twice-differentiable function $r: \real \to \real$ with $r(0) = 0$. We assume that $R 
> 1$, that $r'(x_1) \ge 1 + \epsilon$ for some $\epsilon > 0$, and that $r''(x_1)$ is 
bounded for all $x_1 \in \real$. (This assumption allows 
us to use the $\mathcal L_\infty$ norm for simplified analysis; weighted norms can be 
used to generalize the parameter ranges.) The state variables are $x_1 \in \real$ 
(voltage across the nonlinear element) and $x_2 \in \real$ (current through the 
inductor), and the control input $u \in \real$ is the voltage across the source. The 
dynamics are
\[
	\dot x_1 = \frac{1}{C} \left(-r(x_1) + x_2\right), \quad
	\dot x_2 = \frac{1}{L} \left(-x_1 - R x_2 + u\right)
\]
from an initial condition $x(0) = \vect 0_2$. Our objective is to minimize the cost
\[
	J[u] = \int_0^T \underbrace{\frac 1 2 u^2(t)}_{\phi}\,dt + 
	\underbrace{\frac{\gamma}{2}
	\norm{x(T) - x^*}{2}^2}_{\psi}
\]
for some terminal cost weight $\gamma > 0$, where $x^* \in \real^2$ is an arbitrary 
target state, and the space of permissible controls is $U = [-u_{\rm max}, u_{\rm max}]$ 
for some $u_{\rm max} > 0$. Note that $\vect 0_2$ is an equilibrium point of the 
unforced 
dynamics. 

\subsection{Examining the Assumptions}
This optimal control problem satisfies Assumptions 
\ref{ass:contraction}--\ref{ass:opt-lip}, as we demonstrate in the following paragraphs. 

\paragraph*{Assumption \ref{ass:contraction}} The dynamics are strongly infinitesimally 
contracting with respect to the $\mathcal L_\infty$ norm:
\begin{align*}
	\osL(f) &= \sup_{x \in \real^2} \mu_\infty(D_x f(x)) \\
	&= \sup_{x_1 \in \real} \max\left\{
		\frac{1 -r'(x_1)}{C}, \; \frac{1 - R}{L}
	\right\} \\
	&= \max\left\{
		\frac{1 - d_{\rm min}}{C}, \; \frac{1 - R}{L}
	\right\} \triangleq -c < 0
\end{align*}
where $d_{\rm min} = \inf_{x_1 \in \real} r'(x_1) > 1$. Thus Assumption 
\ref{ass:contraction} is satisfied with contraction rate $c$.

\paragraph*{Assumption \ref{ass:lip-1}} Given two inputs 
$u, \bar u \in \real$, $\norm{f(x, u) - f(x, \bar u)}{\infty} = L^{-1} |u - \bar 
u|$ for all $x \in \real^2$, so Assumption \ref{ass:lip-1} is satisfied with $\ell_{f, 
u} 
= L^{-1}$. 

\paragraph*{Reachability Analysis} Before we examine Assumption \ref{ass:lip-2}, it is 
useful to bound the set of states that are reachable within time $T$. Lemma 
\ref{lem:iss} 
allows us to compare $x(t)$ with the trajectory at the origin corresponding to zero 
input:
\[
	\norm{x(t)}{\infty} \le L^{-1} \int_0^t e^{-c(t - \tau)} |u(\tau)|~d\tau \le 
	\frac{u_{\rm max}(1 - e^{-c t})}{c L}  
\]
In particular, $x(T)$ belongs to a $\mathcal L_\infty $ ball centered about the origin, 
with radius $u_{\rm max} L^{-1} \kappa$, where $\kappa$ is defined in \eqref{eq:kappa}.

\paragraph*{Assumption \ref{ass:lip-2}} 
Since the Jacobian matrix $D_x f(x, u)$ has no dependence on $u$, we have $\ell_{f_x, u} 
= 0$. To evaluate $\ell_{f_x, x}$, note that
\begin{align*}
	D_x f(x, u)^\top \lambda &= \begin{bmatrix} -C^{-1} r'(x_1) & -L^{-1} \\  C^{-1} & -R 
	L^{-1} \end{bmatrix} \begin{bmatrix} \lambda_1 \\ \lambda_2 \end{bmatrix}  
\end{align*}
Then for any $x, \bar x \in \real^2$, 
\begin{align*}
	\norm{D_x f(x, u)^\top \lambda - D_x f(\bar x, u)^\top \lambda}{1} &= C^{-1} 
	|\lambda_1| |r'(x_1) - r'(x_2)| \\
	&\le C^{-1} \eta |\lambda_1| |x_1 - x_2|
\end{align*}
where we define $\eta = \sup_{x \in \mathbb R} |r''(x)|$.
(Note we evaluate the $\mathcal L_1$ norm, which is dual to the $\mathcal L_\infty$ norm 
of the state space.) To bound $|\lambda_1|$, we note that Theorem \ref{thm:dual} implies 
that the time-reversed costate dynamics $\LambBack(t)$ are strongly infinitesimally 
contracting with rate $c$. Furthermore, the origin is a trajectory, so by Lemma 
\ref{lem:iss},
\[
	\norm{\LambBack(t)}{1} \le e^{-c t} \norm{\LambBack(0)}{1} \le 
	\norm{\lambda(T)}{1}, \; \forall t \in [0, T].
\]
Since $\lambda(T) = \psi_x(x(T))$, we can then bound
\[
	\norm{\lambda(t)}{1} \le \norm{\psi_x(x(T))}{1} 
	= \gamma \norm{x(T) - x^*}{1}, \; \forall t \in [0, T].
\]  
Then for all $t \in [0, T]$,
\[
	|\lambda_1(t)| \le \norm{\lambda(t)}{1} \le \gamma \left(
		\norm{x^*}{1} + 2 u_{\rm max} L^{-1} \kappa
	\right),
\]
using the property that $\norm{x(T)}{\infty} \le u_{\rm max} L^{-1} \kappa$. Thus, 
Assumption \ref{ass:lip-2} is satisfied with
\[
	\ell_{f_x, x} = \frac{\gamma \eta}{C} \left(\norm{x^*}{1} + 2 u_{\rm max} L^{-1} 
	\kappa \right).
\]

\paragraph*{Assumption \ref{ass:input-lip}}

Since the running cost $\phi(u) = u^2$ has no dependence on $x$, we have 
$\ell_{\phi_x,x} 
= 0$ and $\ell_{\phi_x, u} = 0$. Furthermore, for any $x, \bar x \in \real^2$, 
$\norm{\psi_x(x) - \psi_x(\bar x)}{1} = \gamma \norm{x - \bar x}{1}$, so Assumption 
\ref{ass:input-lip} is satisfied with $\ell_{\psi_x, x} = \gamma$.

\paragraph*{Assumption \ref{ass:opt-lip}}

The Hamiltonian can be written
\[
	H(x, \lambda, u) = \frac{1}{2} u^2 + \frac{\lambda_2}{L} u + b(x, \lambda)
\]
for a constant offset $b(x, \lambda)$. Minimizing the Hamiltonian over $u \in [-u_{\rm 
max}, u_{\rm max}]$ leads to the unique minimizer
\[
	h(\lambda) = \begin{cases}
		-u_{\rm max}, & L^{-1} \lambda_2 > u_{\rm max} \\
		-L^{-1} \lambda_2, & L^{-1} |\lambda_2| \le u_{\rm max} \\
		u_{\rm max}, & L^{-1} \lambda_2 < -u_{\rm max}
	\end{cases}
\]
The map $h$ is Lipschitz in $\lambda$ with no dependence on $x$, so Assumption 
\ref{ass:opt-lip} is satisfied with $\ell_{h, x} = 0$ and $\ell_{h, \lambda} = L^{-1}$.

\subsection{Convergence of MSA}

Having demonstrated that the optimal control problem satisfies Assumptions 
\ref{ass:contraction}--\ref{ass:opt-lip}, we can state the guarantees of Theorem 
\ref{thm:main}. Substituting in the Lipschitz constants from the previous section into 
\eqref{eq:b1}--\eqref{eq:b2}, we obtain
\[
	b_1 = \frac{\gamma}{L^2}, \quad	b_2 = \frac{\gamma \eta}{C L^2} \left( 
	\norm{x^*}{1} + 2 u_{\rm max} L^{-1} \kappa \right) 
\]
By Theorem \ref{thm:main}, convergence is guaranteed when 
\[
	\kappa + \frac{\eta \norm{x^*}{1}}{C} \kappa^2 + \frac{2 \eta u_{\rm max}}{L 
	C} \kappa^3 < \frac{L^2}{\gamma}
\]  

\subsection{Numerical Results}

\begin{figure}
	\centering
	\includegraphics[width=0.9\linewidth]{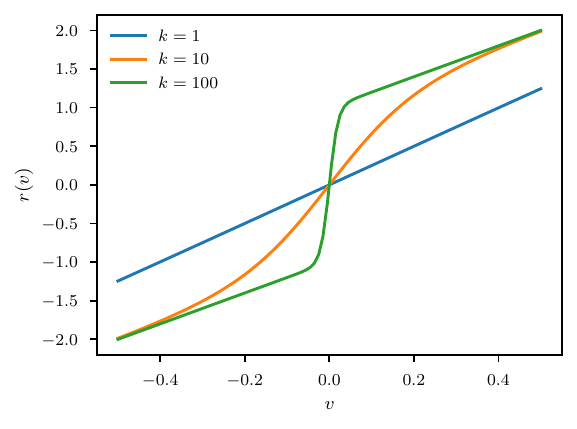}
	\caption{Shape of the $r(v)$ function from \eqref{eq:r}, with $\alpha = \beta = 2$, 
	for various values of $k$.}
	\label{fig:r-nonlin}
\end{figure}

\begin{figure}
	\centering
	\includegraphics[width=0.9\linewidth]{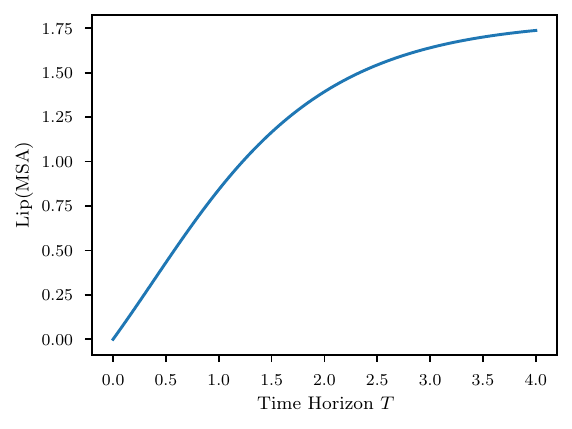}
	\caption{Bounds on the Lipschitz constant of the MSA operator at various time 
	horizons, via Theorem \ref{thm:main}.}
	\label{fig:msa-lip}
\end{figure}

\begin{figure}
	\centering
	\includegraphics[width=0.9\linewidth]{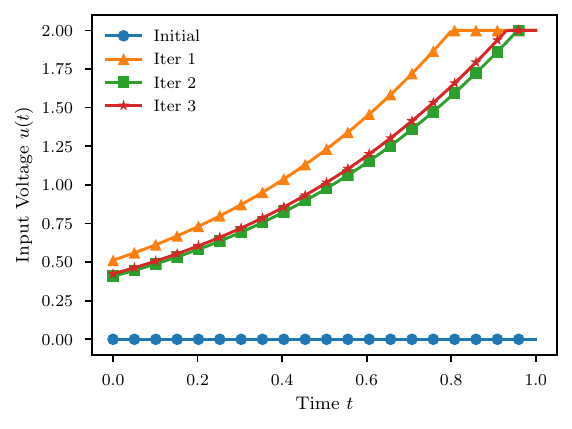}
	\caption{Three successive iterates of the MSA algorithm when $\Lip(\MSA) \le 0.85$, 
	starting from an initial guess $u^{(0)}(t) = 0$.}
	\label{fig:inputs-converge}
\end{figure}

Consider a nonlinearity of the form
\begin{equation} \label{eq:r}
	r(v) = \alpha v + \beta \left( \frac{1}{1 - e^{-k v}} - \frac 1 2\right)
\end{equation}
so that $d_{\rm min} \triangleq \inf_{v \in \real} r'(v) = \alpha$ and $\eta \triangleq 
\sup_{v \in \real} |r''(v)| = \beta k^2 / (6 \sqrt 3)$. Figure \ref{fig:r-nonlin} 
illustrates this function for various values of the shape parameter $k$. We select 
$\alpha = \beta = 2$ and $k = 4$. Furthermore, we select a target state $x^* = (0.1, 
r(0.1))$, with $u_{\rm max} = 2$ and terminal cost weight $\gamma = 100$, with model 
parameters $R = 20$, $L = 11$, and $C = 1$. With these parameters, the contraction rate 
is $c = 1$, and the upper bound on $\Lip(\MSA)$ from Theorem \ref{thm:main} is 
plotted in Figure \ref{fig:msa-lip}. We select a time horizon of $T = 1$, where the 
bound $\Lip(\MSA) \le 0.85$ is guaranteed.

In order to implement the MSA algorithm, we use \verb|solve_ivp| from the 
SciPy package to integrate the state and costate dynamics. This function implements the 
explicit Runge-Kutta method RK45, and it approximates the solution as a continuous 
function using quartic interpolation. Starting with an initial guess of $u^{(0)} = 0$, 
the MSA algorithm quickly converges, with the iterates $u^{(i)}$ for the first three 
iterations $i \in \{1, 2, 3\}$ depicted in Figure \ref{fig:inputs-converge}. Note the 
rapid decay of the $\mathcal L_\infty$ distance between each successive iterate.

\section{Conclusions}

In this letter, we have examined an indirect method for the optimal control of strongly 
contracting systems. We have observed that the time-reversed adjoints of such systems are 
also contracting with the same rate, with respect to the dual norm, leading to useful 
bounds on the costate trajectory from PMP. Based on this observation, we bounded the 
Lipschitz constant of each iteration of MSA, demonstrating that the iteration is actually 
a contraction mapping for sufficiently strongly contractive systems or for sufficiently 
short horizons. In these cases, MSA is guaranteed to converge to a unique control that 
satisfies PMP. With an additional assumption on pointwise uniqueness of the minimizer of 
the Hamiltonian, we showed that this control is indeed the unique optimal control.
 
The main approach of this paper, namely using ISS properties of the adjoint to bound the 
Lipschitz constant of a $\mathcal U \to \mathcal U$ operator, is quite general and could 
be applied to many other indirect methods in optimal control. Several variants of MSA, 
both older \cite{SKM:66, FLC-AAL:82} and newer \cite{QL-SH:18, QL-LC-CT-WE:18}, could be 
studied with this type of analysis in future work, possibly with more general 
convergence criteria. Another practical future direction would be the study of 
discretized implementation of the forward and backward integration steps, as in 
\cite{MMcA-LM-WH:12}. Of course, one could also analyze indirect methods for extensions 
of the optimal control problem, such as constraints on the terminal state or in the 
infinite time horizon. An additional interesting direction would be the application of 
convergence guarantees to model predictive control of contractive nonlinear systems.

\appendices
\section{Proofs}

\subsection{Proof of Lemma \ref{lem:dual-lip}}

\rev{
For a general matrix $A \in \real^{n \times n}$, from the definitions of induced 
norms and dual norms we have
\[
	\norm{A}{\star} = \sup_{\norm{y}{\star} = 1} \sup_{\norm{z} \le 1} z^\top A y 
\]
Swapping the order of the suprema and applying the definition of the dual 
norm once again yields
\[
	\norm{A}{\star} = \sup_{\norm{z} \le 1} \sup_{\norm{y}{\star} = 1} y^\top A^\top z
	= \sup_{\norm{z} \le 1} \norm{A^\top z}{\star\star}
\]
But $\real^n$ with any norm is a reflexive Banach space, so $\norm{\cdot}_{\star\star} = 
\norm{\cdot}$, and thus $\norm{A}{\star} = \norm{A^\top}$. We use this fact to prove both 
statements. Because $g$ is continuously differentiable and $D_x g(x) = D_x f(x)^\top$,
\[
	\Lip_{\norm{\cdot}_\star}(g) = \sup_{x \in X} \norm{D_x g(x)}{\star}
	= \sup_{x \in X} \norm{D_x f(x)} = \Lip_{\norm{\cdot}}(f)
\]
and
\begin{align*}
	\osL_{\norm{\cdot}{\star}}(g) &= \sup_{x \in X} \lim_{\alpha \to 0^+} \frac{\norm{I_n 
	+ \alpha D_x g(x)}{\star} - 1}{\alpha} \\
	&= \sup_{x \in X} \lim_{\alpha \to 0^+} \frac{\norm{I_n 
	+ \alpha D_x f(x)} - 1}{\alpha} = \osL_{\norm{\cdot}}(f)
\end{align*}
}  

\subsection{Proof of Lemma \ref{lem:bounded}}

\rev{
Let $\bar x(t)$ be the trajectory of \eqref{eq:sys} corresponding to input $\bar u(t) = 
\vect 0_k$. Since \eqref{eq:sys} is strongly infinitesimally contracting, we can use 
Lemma \ref{lem:iss} to compare a trajectory $x(t)$ with $\bar x(t)$:
\[
	\norm{x(t) - \bar x(t)} \le \frac{\ell_{f, u}}{c} (1 - e^{-cT}) 
	\!\sup_{\tau \in [0, T]}\! \norm{u(\tau)}{U}
\]
for all $t \in [0, T]$. Since $U$ and $\bar x(t)$ are bounded, $x(t)$ is bounded as 
well. Similarly, the time-reversed costate dynamics \eqref{eq:adjoint} have an 
equilibrium point at the origin when $v(t) = \vect 0_n$, regardless of $x(t)$ and 
$u(t)$, and (due to Theorem \ref{thm:dual}) they are strongly 
contracting with rate $c > 0$. Again, we can use Lemma \ref{lem:iss} to compare 
$\LambBack(t)$ with the trajectory at the origin:
\[
	\norm{\LambBack(t)}{\star} \le \norm{\LambBack(0)}{\star} + 
	\frac{1}{c}(1 - e^{-cT}) 
	\sup_{\tau \in [0, T]} \norm{v(\tau)}{\star}
\] 
for all $t \in [0, T]$. Since $V$ is bounded, $\LambBack(t)$ is confined to a ball 
$\Lambda$ about the origin.
}

\subsection{Proof of Theorem \ref{thm:costate-iss}}

As in Theorem \ref{thm:dual}, let $\LambBack(t) = \lambda(T - t)$, so that
\[
	\frac{d \LambBack(t)}{dt} = D_x f(T - t, x(T - t), u(T - t))^\top \LambBack(t) - v(t 
	- T)
\]
At any fixed $t$, the $\LambBack$ vector field has the 
Jacobian matrix $D_x f(T - t, x(T - t), u(T - t))^\top$, which is transpose the 
Jacobian matrix of $f(T - t, \cdot, u(T - t))$. By Assumption \ref{ass:contraction}, 
$\osL(f(T - t, \cdot, u(T - t))) \le -c$, so Lemma \ref{lem:dual-lip} implies that the 
$\LambBack$ vector field is also one-sided Lipschitz with constant $c$, with respect to 
$\norm{\cdot}_\star$. Then we apply Lemma \ref{lem:iss} to bound $\norm{\LambBack(t) - 
\bar \LambBack(t)}_\star$ with respect to the inputs $u(t)$, $x(t)$, and $v(t)$, 
resulting in the following bound on $\norm{\lambda(t) - \bar \lambda(t)}_\star$:
\begin{align*}
	&\norm{\lambda(t) - \bar \lambda(t)}_\star \le e^{-c (T - t)} \norm{\lambda(T) - \bar 
	\lambda(T)}_\star \\
	&\qquad\qquad + \ell_{f_x, u} \int_{t}^{T} e^{-c(\tau - t)} \norm{u(\tau) - 
	\bar u(\tau)}_U ~d\tau \\
	&\qquad\qquad + \ell_{f_x, x} \int_t^{T} e^{-c(\tau - t)} \norm{x(\tau) - 
	\bar x(\tau)}_X~d\tau \\
	&\qquad\qquad + \int_t^{T} e^{-c(\tau - t)} \norm{v(\tau) - \bar v(\tau)}_\star~d\tau
\end{align*}
We apply Lemma \ref{lem:iss} once more to remove explicit dependence on $x$, via the bound
\begin{multline*} 
	\int_t^T e^{-c(\tau - t)} \norm{x(\tau) - \bar x(\tau)}_X~d\tau \\
	\le \ell_{f, u} \!\int_t^T\!\! \!\int_{0}^\tau\!\! e^{-c(\tau - t)} 
	e^{-c(\tau - \tau')} \norm{u(\tau') - \bar u(\tau')}_U~d\tau d\tau'  
\end{multline*}
We then swap the order of integration:
\begin{align*}
	\int_t^T \int_{0}^\tau & e^{-c(\tau - t)} e^{-c(\tau - \tau')} \norm{u(\tau') - \bar 
	u(\tau')}_U d\tau' d\tau \\
	&= \!\int_0^t\! \!\int_t^T\! e^{-c(\tau - t)} e^{-c(\tau - \tau')} \norm{u(\tau') - 
	\bar u(\tau')}_U d\tau d\tau' \\
	&\;\; + \!\int_t^T\! \!\int_{\tau'}^T\! e^{-c(\tau - t)} e^{-c(\tau - \tau')} 
	\norm{u(\tau') - \bar u(\tau')}_U d\tau d\tau' \\
	&= \frac{\sinh(c(T - t))}{c} \!\int_0^t\! e^{-c(T - \tau)} \norm{u(\tau) - \bar u \tau)}_U d\tau \\
	&\;\; + \frac{e^{-c(T - t)}}{c} \!\int_t^T\! \sinh(c(T - \tau))  \norm{u(\tau) - \bar u \tau)}_U d\tau
\end{align*} 

\subsection{Proof of Corollary \ref{cor:sup-bound}}

The first three terms are obvious upper bounds on the first three terms in 
\eqref{eq:costate-iss}, and
\begin{align*}
	&\frac{\sinh(c(T \!-\! t))}{c} \!\!\! \int_0^t \!\!\! e^{-c(T \!-\! \tau)}\;d\tau 
	+ \frac{e^{-c(T \!-\! t)}}{c} \!\!\! \int_t^T \!\!\! \sinh(c(T \!-\! 
	\tau))\;d\tau \\
	&= \!\int_t^T\! \!\int_0^\tau\! e^{-c(\tau - t)} e^{-c(\tau - \tau')}~d\tau' d\tau 
	\le \kappa \!\int_t^T\! e^{-c(\tau - t)} ~d\tau \le \kappa^2 
\end{align*}  

\subsection{Proof of Theorem \ref{thm:main}}

Let $u, \bar u \in \mathcal U$, and let $x, \bar x$ and $\lambda, \bar \lambda$ be 
the corresponding state and costate trajectories. Then for all $t \in [0, T]$,
\begin{align}
	\begin{split} \label{eq:fbs-bound}
	&\norm{\MSA(u)(t) - \MSA(\bar u)(t)}_U \\
	&\qquad = \norm{h(t, x(t), \lambda(t)) - h(t, \bar x(t), \bar \lambda(t))}_U \\
	&\qquad \le \ell_{h, x} \norm{x(t) - \bar x(t)} + \ell_{h, \lambda} \norm{\lambda(t) - 
	\bar \lambda(t)}_\star
	\end{split}
\end{align}
The costate dynamics are \eqref{eq:adjoint} with $v(t) = -\phi_x(t, x(t), u(t))$, which 
is bounded in $(\real^n, \norm{\cdot}{\star})$ by the boundedness of $x(t)$ and $u(t)$ 
and the Lipschitz continuity of $\phi_x$. By Corollary \ref{cor:sup-bound},
\begin{align*}
	&\norm{\lambda(t) - \bar \lambda(t)}_\star \le \norm{\lambda(T) - \bar \lambda(T)}_\star \\
	&\qquad + \kappa \sup_{t \in [0, T]} \norm{\phi_x(t, x(t), u(t)) - \phi_x(t, \bar 
	x(t), \bar u(t))}_\star \\
	&\qquad + \left( \ell_{f_x, u} \kappa + \ell_{f_x, x} \ell_{f, u} \kappa^2 \right)
	\sup_{t \in [0, T]} \norm{u(t) - \bar u(t)}_U,
\end{align*}
where
\begin{align*}
  &\norm{\lambda(T) - \bar \lambda(T)}_\star = \norm{\psi_x(x(T)) - \psi_x(\bar 
    x(T))}_\star \\
  &\le \ell_{\psi_x, x} \norm{x(T) - \bar x(T)}
  \le \ell_{\psi_x, x} \sup_{t \in [0, T]} \norm{x(t) - \bar x(t)}
\end{align*}
and
\begin{align*}
  & \norm{\phi_x(t, x(t), u(t)) - \phi_x(t, \bar x(t), \bar u(t))}_\star \\
  &\qquad \le \ell_{\phi_x, x} \norm{x(t) - \bar x(t)} + \ell_{\phi_x, u} \norm{u(t) - 
  \bar u(t)}_U \\
  &\qquad \le \ell_{\phi_x, x} \!\!\sup_{t \in [0, T]}\! \norm{x(t) - \bar x(t)} \! + \! 
  \ell_{\phi_x, u} \!\!\sup_{t \in [0, T]}\! \norm{u(t) - \bar u(t)}_U
\end{align*}	
As a consequence of Lemma \ref{lem:iss}, $\sup_{t \in [0, T]} \norm{x(t) - \bar x(t)} \le 
\ell_{f, u} \kappa \sup_{t \in [0, T]} \norm{u(t) - \bar u(t)}_U$, so we simplify
\begin{align*}
  & \norm{\lambda(t) - \bar \lambda(t)}_\star \le
  \ell_{\psi_x, x} \ell_{f, u} \kappa \sup_{t \in [0, T]} \norm{u(t) - \bar u(t)}_U \\
  &\qquad + \kappa^2 \ell_{\phi_x, x} \ell_{f, u} \sup_{t \in [0, T]} \norm{u(t) - \bar 
  u(t)}_U \\
  &\qquad + \kappa \ell_{\phi_x, u} \sup_{t \in [0, T]} \norm{u(t) - \bar u(t)}_U \\
  &\qquad + (\ell_{f_x, u} \kappa + \ell_{f_x, x} \ell_{f, u} \kappa^2) \sup_{t \in [0, 
      T]} \norm{u(t) - \bar u(t)}_U.
\end{align*}
Substituting the state and costate difference bounds into \eqref{eq:fbs-bound} completes 
the proof of statement \ref{item:lip}. Then statement \ref{item:banach} is a standard 
consequence of the Banach fixed point theorem.

\subsection{Proof of Corollary \ref{corollary}}

We first establish that the fixed points of the MSA operator are precisely the 
controls that satisfy PMP. One direction is obvious: $u^* = \MSA(u^*)$ implies that 
$u^*$ satisfies PMP. Now suppose that $u^*$ satisfies PMP, and let $x^*, \lambda^*$ 
be the corresponding state and costate trajectories. Then $u^*(t) \in \argmin_{u \in 
U} H(t, x^*(t), \lambda^*(t), u)$ for all $t \in [0, T]$, so the assumption that the 
Hamiltonian has a unique minimizer implies that $u^*(t) = h(t, x^*(t), \lambda^*(t))$ for 
all $t \in [0, T]$, and thus $u^* = \MSA(u^*)$.

We then establish that the MSA iteration converges to a unique fixed point $\hat 
u$. For $T$ sufficiently small or $c$ sufficiently large, $\kappa$ is sufficiently 
small that $\Lip(\MSA) \le b_1 \kappa + b_2 \kappa^2 < 1$, by Theorem \ref{thm:main}. 
Then the Banach fixed point theorem establishes that a unique fixed point $\hat u$ 
exists, and that the iteration from any initial guess converges to $\hat u$. 

Since an optimal control $u^*$ exists, it is a fixed 
point of $\MSA$, and the fixed point of $\MSA$ is unique. Furthermore, if a control $u^*$ 
satisfies PMP, then it is a fixed point of $\MSA$, and hence is equal to the optimal 
control.
	
\bullobib
 
\end{document}